
\input amssym.def
\input amssym
\input psfig
\magnification=1100
\baselineskip = 0.25truein
\lineskiplimit = 0.01truein
\lineskip = 0.01truein
\vsize = 8.5truein
\voffset = 0.2truein
\parskip = 0.10truein
\parindent = 0.3truein
\settabs 12 \columns
\hsize = 5.4truein
\hoffset = 0.4truein
\font\ninerm=cmr9

\setbox\strutbox=\hbox{%
\vrule height .708\baselineskip
depth .292\baselineskip
width 0pt}
\font\caps=cmcsc10
\font\smallmaths=cmmi9
\font\normalmaths=cmmi10
\font\bigtenrm=cmr10 at 14pt

\def\sqr#1#2{{\vcenter{\vbox{\hrule height.#2pt
\hbox{\vrule width.#2pt height#1pt \kern#1pt
\vrule width.#2pt}
\hrule height.#2pt}}}}
\def\square{\mathchoice\sqr46\sqr46\sqr{3.1}6\sqr{2.3}4}

\centerline{\bigtenrm SOME 3-MANIFOLDS AND 3-ORBIFOLDS}
\centerline{\bigtenrm WITH LARGE FUNDAMENTAL GROUP}
\tenrm
\vskip 14pt
\centerline{MARC LACKENBY}
\vskip 18pt

\centerline{\caps Abstract}
\vskip 6pt

\ninerm
\textfont1 = \smallmaths
\baselineskip = 0.22truein
We provide two new proofs of a theorem of Cooper, Long
and Reid which asserts that, apart from an explicit finite
list of exceptional manifolds, any compact orientable irreducible
3-manifold with non-empty boundary has large fundamental group.
The first proof is direct
and topological; the second is group-theoretic. These techniques are 
then applied to prove a string of results about (possibly closed) 3-orbifolds,
which culminate in the following theorem. If $K$ is a knot
in a compact orientable 3-manifold $M$, such that the complement
of $K$ admits a finite-volume hyperbolic structure,
then the orbifold obtained by assigning a singularity of
order $n$ along $K$ has large fundamental group, for infinitely
many positive integers $n$. We also obtain information
about this set of values of $n$. When $M$ is the 3-sphere, this
has implications for the cyclic branched covers over the knot. In this
case, we may also weaken the hypothesis that the complement
of $K$ is hyperbolic to the assumption that $K$ is non-trivial.

\tenrm
\textfont1 = \normalmaths
\baselineskip = 0.25truein
\vskip 18pt
\centerline {\caps 1. Introduction}
\vskip 6pt

One of the key unresolved problems in 3-manifold theory
is the Virtually Haken Conjecture. This asserts that
if $M$ is a closed orientable hyperbolic 3-manifold, 
then $M$ has a finite-sheeted
cover containing a properly embedded orientable incompressible
surface (other than a 2-sphere). There are a range of stronger
forms of this conjecture, possibly the strongest of
which proposes that the fundamental group of
$M$ is `large'. This means that it has a finite
index subgroup that admits a surjective homomorphism
onto a non-abelian free group. The covering space of $M$ corresponding
to such a finite index subgroup contains the required
incompressible surface. But large groups have many
other nice properties. For example, they have finite
index subgroups with arbitrarily large first Betti number.
In this paper, we will exhibit several classes of 3-manifolds
and 3-orbifolds
with large fundamental group. We start with a rapid
and surprisingly elementary proof of the following theorem of 
Cooper, Long and Reid. 

\noindent {\bf Theorem 2.1.} [5] {\sl Let $M$ be a compact orientable
irreducible 3-manifold with non-empty boundary. Then, either
$M$ is an $I$-bundle over a surface with non-negative
Euler characteristic or $\pi_1(M)$ is large.}

The original proof by Cooper, Long and Reid relied
on some delicate and complicated 3-dimensional
techniques. However, their aim was somewhat different
from ours. They first showed that a compact orientable 
hyperbolic 3-manifold $M$ with non-empty boundary
(other than an $I$-bundle over a surface)
has a finite-sheeted cover that
contains a closed embedded orientable incompressible surface
(other than a 2-sphere) that is not boundary parallel. Further
work was then required to ensure that this
surface is non-separating, and yet more work
was necessary to find two such surfaces whose
union is non-separating. This then implies that
$\pi_1(M)$ is large. Establishing these
results were hard work, particularly when $\partial M$
does not consist entirely of tori. However, if one is
content solely with proving 
that $\pi_1(M)$ is large, without establishing the existence
of a closed essential surface in some finite cover, 
then most of the difficult 3-dimensional arguments
can be sidestepped, as we shall see.

It is interesting to note that, in fact, one can almost entirely
dispense with the topology and prove Theorem 2.1
using primarily group-theoretic techniques (at least in
the main case, when $M$ is hyperbolic). We supply
such a proof, which uses a recent theorem of the author, that
gives a criterion for a finitely presented group to be large, 
in terms of the behaviour of its finite
index subgroups [9]. However, there is still
topology lurking in the background, as this
largeness criterion was proved using topological
methods.

One of the main reasons why the topological proof of Theorem 2.1
is useful is that it naturally extends to certain
3-orbifolds. We pursue this line of investigation in 
\S3, and give a brief description of these results below.
These provide new classes of 3-orbifolds (and hence 3-manifolds) 
that have large fundamental group. 

Throughout this paper, an orbifold is allowed
to have empty singular locus, and hence be a manifold.
If $O$ is a 3-orbifold and $L$ is a link in $O$ disjoint
from the singular locus and $n$ is a positive integer, 
then we denote by $O(L,n)$ the orbifold
obtained from $O$ by adding singularities along $L$ of order $n$.
We will prove a sequence of results about 3-orbifolds, 
which lead to the following.

\noindent {\bf Theorem 3.6.} {\sl Let $O$ be a
compact orientable 3-orbifold (with possibly empty singular locus), and
let $K$ be a knot in $O$, disjoint from its singular
locus, such that $O - K$ has a finite volume
hyperbolic structure. Then, for infinitely many
values of $n$, $\pi_1(O(K,n))$ is large.}

Setting $O$ to be the 3-sphere in Theorem 3.6, this applies
to the much studied case of cyclic branched covers
over hyperbolic classical knots. In fact, by use of the
Orbifold Theorem and applying results from \S3, we can obtain
the following information about branched covers over any non-trivial
knot in the 3-sphere.

\noindent {\bf Theorem 3.7.} {\sl Let $K$ be a
non-trivial knot in the 3-sphere, and let $m$ be any integer
more than two. Then, for all sufficiently large $n$, the $mn$-fold
cyclic cover of $S^3$ branched over $K$ has
large fundamental group.}

Some parts of this paper present new proofs of 
known results; other bits give new theorems.
However, the outstanding paper of Cooper,
Long and Reid [5] has exerted a strong influence
throughout.

Another simple proof of Theorem 2.1 has appeared
recently, due to Button [3]. He showed that it can
be deduced quite quickly from Howie's criterion [7]
for a group to be large. Ratcliffe also established
largeness in  the case where $M$ has a boundary
component with genus at least two in [14], providing
a very quick proof based on a theorem of Baumslag
and Pride [1].

\vskip 18pt
\centerline{\caps 2. Bounded 3-manifolds}
\vskip 6pt

Our goal in this section is to provide two new proofs of the following theorem.

\noindent {\bf Theorem 2.1.} [5] {\sl Let $M$ be a compact orientable
irreducible 3-manifold with non-empty boundary. Then, either
$M$ is an $I$-bundle over a surface with non-negative
Euler characteristic or $\pi_1(M)$ is large.}

The principle reason why compact orientable 3-manifolds with non-empty
boundary are more tractable than closed manifolds
is the following result. This gives a lower bound on the 
rank of their cohomology in terms of the genus of their boundary.
This result is a well known consequence of Poincar\'e
duality.

\vfill\eject
\noindent {\bf Proposition 2.2.} {\sl Let $M$ be a compact
orientable 3-manifold, and let $i \colon P \rightarrow M$ be 
the inclusion of a compact (possibly empty) subsurface
of $\partial M$. Then
$${\rm rank}({\rm ker}(i^\ast \colon H^1(M) \rightarrow H^1(P))) \geq 
{\textstyle{1 \over 2}} b_1(\partial M)
- b_1(P).$$}

To establish largeness, we will use the following 
well known lemma.

\noindent {\bf Lemma 2.3.} {\sl Let $M$ be a compact
3-manifold. Suppose that $M$ contains two disjoint, transversely oriented,
properly embedded surfaces whose union is non-separating. Then
$\pi_1(M)$ admits a surjective homomorphism onto ${\Bbb Z} \ast {\Bbb Z}$.}

\noindent {\sl Proof.} Let $S_1$ and $S_2$ be the two surfaces.
There is the following collapsing map
$f \colon M \rightarrow S^1 \vee S^1$. The restriction of $f$
to a regular neighbourhood $N(S_i)$ is the composition of the
homeomorphism $N(S_i) \rightarrow S_i \times I$
with projection onto the $I$ factor, followed
by the quotient map from $I$ to the circle
that glues the ends of the interval together,
composed with the inclusion into the $i^{\rm th}$
circle of $S^1 \vee S^1$.
The map sends the remainder of $M$ to the central
vertex of $S^1 \vee S^1$. Fix a basepoint
in $M$ disjoint from $N(S_1) \cup N(S_2)$. It is clear that
$f_\ast \colon \pi_1(M) \rightarrow {\Bbb Z} \ast {\Bbb Z}$
is a surjection, because any element of
${\Bbb Z} \ast {\Bbb Z}$ may be realised by a based
loop in $M$. $\square$

\noindent {\sl Proof of Theorem 2.1.} 
First suppose that
$\partial M$ contains a 2-sphere. Then, by irreducibility,
$M$ is a 3-ball, which is an $I$-bundle over a disc,
verifying the theorem in this case. Thus, we may
assume that each component of $\partial M$ has
genus at least one.

Consider first the case where
each component of $\partial M$ is a torus.
A standard argument then allows us to
assume that $M$ is hyperbolic. This argument
can be found in Cooper, Long and Reid's paper [5],
but we repeat it here for the sake of completeness.
When $M$ is not hyperbolic, Thurston's geometrisation
theorem [13] implies that $M$ is either Seifert
fibred or contains an essential embedded torus.
In the former case, the argument divides according
to whether the base orbifold of the Seifert fibration has positive, zero
or negative Euler characteristic. When it
is positive, the manifold is a solid torus, which
is an $I$-bundle over an annulus. When it is
zero, the manifold again admits some $I$-bundle
structure over a torus or Klein bottle. When
the base orbifold has negative Euler characteristic, it has a finite-sheeted
cover which is an orientable surface also with
negative Euler characteristic.
This induces a finite covering $\tilde M \rightarrow M$.
The Seifert fibration on $\tilde M$ induces a surjective homomorphism
from $\pi_1(\tilde M)$ onto the fundamental group of
this surface, which then admits a surjective
homomorphism onto a free non-abelian group. Hence,
when $M$ is Seifert fibred, it satisfies the
conclusion of the theorem. 

If $M$ contains an essential embedded torus, then it is a theorem
of Kojima [8] and Luecke [12] (see also the work of Niblo 
and Long in [10] and [11]) that either $\pi_1(M)$ is large
or $M$ is finitely covered by a torus bundle over the
circle, the torus lifting to a fibre. The latter case
cannot arise since $M$ has non-empty boundary.

We may therefore assume (when $\partial M$ consists entirely of tori)
that the interior of
$M$ is a finite-volume hyperbolic 3-manifold. A lemma of Cooper, Long and
Reid (Lemma 2.1 of [5]) asserts that, by passing to a finite-sheeted
covering space if necessary, we may assume that $M$ has at least three 
boundary components.

All the above follows the argument of Cooper, Long and Reid,
but here our proofs diverge.
We have reached the stage where {\sl either}
$\partial M$ consists of tori and there are at least
three of these, {\sl or} $\partial M$ has a component
with genus at least two. In the former case, set
$P$ to be one these tori; otherwise let $P$ be
the empty set. Then Proposition 2.2 
gives that the kernel of $i^\ast \colon H^1(M) \rightarrow H^1(P)$
has rank at least ${1 \over 2} b_1(\partial M) - b_1(P)$,
which is positive. Let $\alpha$ be
a non-trivial primitive element in the
kernel of $i^\ast \colon H^1(M) \rightarrow H^1(P)$. Let $S$ be a properly
embedded oriented surface in $M$, disjoint from $P$,
dual to $\alpha$. We may assume that $S$ intersects each toral component
of $\partial M$ in a (possibly empty) collection of
coherently oriented essential curves.
Now, $\alpha$ induces a surjective
homomorphism $\pi_1(M) \rightarrow {\Bbb Z}$.
Composing this with the homomorphism
${\Bbb Z} \rightarrow {\Bbb Z}/n{\Bbb Z}$ that
reduces modulo $n$, we obtain a homomorphism
$\pi_1(M) \rightarrow {\Bbb Z}/n{\Bbb Z}$.
Let $M_n$ be the corresponding
$n$-fold cyclic cover of $M$.
This contains $n$ disjoint copies of $S$ which
can be labelled with the integers modulo $n$.
Let $F_n$ be the union of the surface with
label $0$ and the surface with label $\lfloor n/2 \rfloor$.
This is a separating surface, dividing $M_n$
into two 3-dimensional submanifolds, which we will call $A^1_n$ and $A^2_n$.

When $\partial M$ consists only of tori, $A^1_n$
and $A^2_n$ each contain at least $\lfloor n/2 \rfloor$
copies of $P$. Hence, for $j = 1$ and $2$,
$b_1(\partial A^j_n) \rightarrow \infty$ 
as $n \rightarrow \infty$.
When $\partial M$ contains a non-toral
component, the same is true, since $\partial A^j_n$
consists of $F_n$ and at least $\lfloor n/2 \rfloor$
copies of $\partial M$ cut along $S$, with their boundary
components glued in pairs.
However, in both cases, $b_1(F_n)$
remains independent of $n$. Therefore,
when $n$ is sufficiently large,
${1 \over 2} b_1(\partial A^j_n) > b_1(F_n)$
for both $j = 1$ and $2$.
Let us fix an integer $n$ where these inequalities hold. 
Proposition 2.2 then gives that the kernel
of $i^\ast \colon H^1(A^j_n) \rightarrow
H^1(F_n)$ is non-trivial. Let $W^j_n$ be a 
connected oriented surface, properly embedded in $A_n^j$,
disjoint from $F_n$, dual to a non-trivial primitive class in this
kernel. (See Figure 1 for the case where $\partial M$ is a union of tori.) 
Then $W_n^1$ and $W_n^2$ are disjoint
oriented surfaces the union of which is non-separating
in $M_n$. By Lemma 2.3, $\pi_1(M)$ is large. $\square$

\vskip 18pt
\centerline{\psfig{figure=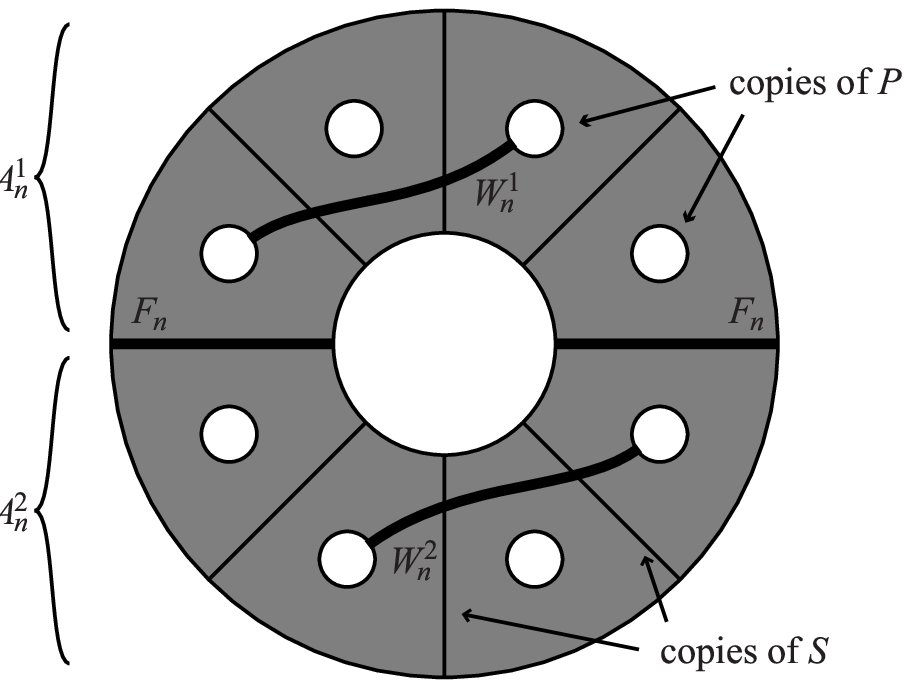,width=3in}}
\vskip 18pt
\centerline{Figure 1.}

We now provide an alternative group-theoretic proof of
Theorem 2.1 that uses the following recent theorem of the author 
(Theorem 1.2 of [9]).

\noindent {\bf Theorem 2.4.} {\sl Let $G$ be a
finitely presented group, and suppose that,
for each natural number $n$, there is a
triple $H_n \geq J_n \geq K_n$ of finite index
normal subgroups of $G$ such that
\item{(i)} $H_n/J_n$ is abelian for all $n$;
\item{(ii)} $\lim_{n \rightarrow \infty} 
((\log [H_n : J_n]) / [G:H_n]) = \infty$;
\item{(iii)} $\limsup_n (d(J_n/K_n) / [G:J_n])  > 0$.

\noindent Then $K_n$ admits a surjective homomorphism
onto a free non-abelian group, for infinitely
many $n$.}

Here, $d( \ )$ denotes the rank of a group, which
is the minimal size of a generating set.

\noindent {\sl Alternative proof of Theorem 2.1.} 
As before, we may assume either that $\partial M$
consists entirely of tori and there are at least
three of these, or that $\partial M$ contains
a higher genus component. Let
$M_n$ be as in the previous proof.

We now wish to apply Theorem 2.4. Let $G$ be
$\pi_1(M)$, and set $H_n$ to be $G$, for each $n$. Let $J_n$
be the subgroup of $G$ corresponding to the
cover $M_n$, and set $K_n$ to be $J_n^2$,
the subgroup generated by the squares of elements
in $J_n$. We must check the various conditions 
of Theorem 2.4. Now, $J_n$ is a normal subgroup of
$G$ by construction. As $K_n$ is a characteristic subgroup of
$J_n$, it is therefore also normal in $G$.
Now, $H_n/J_n$ is isomorphic to ${\Bbb Z}/n{\Bbb Z}$.
In particular, it is abelian, verifying (i),
and its order tends to infinity which gives (ii).
Finally, $J_n/K_n$ is isomorphic to 
$H_1(M_n; {\Bbb Z}/2{\Bbb Z})$, which
has rank at least $n$, by Proposition 2.2. Therefore,
$d(J_n/K_n)/[G:J_n] \geq 1$, for each $n$,
which gives (iii). Hence, by Theorem 2.4,
$\pi_1(M)$ is large. $\square$

\vskip 18pt
\centerline {\caps 3. Orbifolds with large fundamental group}
\vskip 6pt

In this section, we show that the techniques
in the topological proof of Theorem 2.1 can be applied
to 3-orbifolds. This is because we dealt in \S2 with
cyclic covers of large degree, and orbifolds
also have have such covers, provided their
singularities has sufficiently large singularity
order.

We start with the following theorem, where the orbifold
contains at least three distinguished components of
its singular locus. All subsequent theorems will
follow from this result.

\noindent {\bf Theorem 3.1.} {\sl Let $O$ be a compact orientable 3-orbifold,
and let $L$ be a link in $O$, disjoint from the
singular locus of $O$, with at least three components. Then, for all
sufficiently large $n$, $\pi_1(O(L,n))$ is large.}

\noindent {\sl Proof.} 
Let $|O|$ denote the underlying manifold of $O$.
The natural map $O \rightarrow |O|$ induces
a surjective homomorphism $\pi_1(O(L,n)) \rightarrow
\pi_1(|O|(L,n))$. Hence, if the latter group
is large, then so is the former. It therefore suffices
to consider the case where $O$ is a manifold $M$, say.

Let $X$
denote the 3-manifold $M - {\rm int}(N(L))$.
Pick a component $L_1$ of $L$. 
Since $L$ has at least three components, Proposition 2.2
implies that $i^\ast \colon H^1(X) \rightarrow H^1(\partial N(L_1))$
has non-trivial kernel. Let $S$ be a properly embedded,
connected, oriented surface in $X$ disjoint from $\partial N(L_1)$
dual to a non-trivial primitive element in this kernel.
We may assume that $S$ intersects each toral component
of $\partial X$ in a (possibly empty) collection of
coherently oriented essential curves.
Let $X_n \rightarrow X$ denote the associated
$n$-fold cyclic cover of $X$. This extends to a cover $O_n
\rightarrow M(L,n)$. Note that there is an inclusion
map $X_n \rightarrow O_n$.
There are $n$ copies of $S$ in $X_n$ which we
may label with the integers modulo $n$.
Let $F_n$ be the union of the surfaces labelled $0$ and $\lfloor n/2 \rfloor$.
Then $F_n$ divides $X_n$ into two manifolds $A^1_n$
and $A^2_n$. Each contains at least $\lfloor n/2 \rfloor$
copies of $\partial N(L_1)$. For $j=1$ and $2$, let
$P_n^j$ be the copy of $F_n$ in $A_n^j$, together with
any adjacent annuli of $\partial A_n^j - F_n$.
Then $b_1(P_n^j)$ is independent of $n$.
So, Proposition 2.2 implies that,
when $n$ is sufficiently large, $i^\ast \colon
H^1(A^j_n) \rightarrow H^1(P_n^j)$ has non-trivial
kernel, for $j = 1$ and $2$.
Let $W^j_n$ be a connected properly embedded oriented surface in $A^j_n$,
disjoint from $P_n^j$, dual to a primitive class in this kernel. 
Then $W^1_n$ and $W^2_n$ are disjoint oriented surfaces, the union of
which is non-separating in $X_n$. Therefore,  by Lemma 2.3, there is
a surjective homomorphism $\pi_1(X_n) \rightarrow
{\Bbb Z} \ast {\Bbb Z}$. Now, the inclusion map
$X_n \rightarrow O_n$ induces a surjective
homomorphism $\pi_1(X_n) \rightarrow \pi_1(O_n)$.
There is an obvious surjective homomorphism
${\Bbb Z} \ast {\Bbb Z} \rightarrow
({\Bbb Z}/n{\Bbb Z}) \ast ({\Bbb Z}/n{\Bbb Z})$
that respects the free factors.
We claim that $\pi_1(X_n) \rightarrow
{\Bbb Z} \ast {\Bbb Z}$ descends to a surjective
homomorphism $\pi_1(O_n) \rightarrow
({\Bbb Z}/n{\Bbb Z}) \ast ({\Bbb Z}/n{\Bbb Z})$.
To prove this, we must show that the kernel of 
$\pi_1(X_n) \rightarrow \pi_1(O_n)$ is sent to the
identity in $({\Bbb Z}/n{\Bbb Z}) \ast ({\Bbb Z}/n{\Bbb Z})$.
Now, this kernel is normally generated by
powers of the meridian curves that encircle
the singular locus of $O_n$, each power being
the order of the relevant singular component.
The meridian of any singular component not adjacent
to $W_n^1$ or $W_n^2$ is sent to the identity
under $\pi_1(X_n) \rightarrow
{\Bbb Z} \ast {\Bbb Z}$. Each component of the
singular locus adjacent to $W_n^1$ or
$W_n^2$ has order $n$, and the $n^{\rm th}$
power of its meridian is sent to the $n^{\rm th}$
power of one of the free generators
of ${\Bbb Z} \ast {\Bbb Z}$, which is in the
kernel of ${\Bbb Z} \ast {\Bbb Z} \rightarrow
({\Bbb Z}/n{\Bbb Z}) \ast ({\Bbb Z}/n{\Bbb Z})$.
Thus, the claim is verified: there is an induced surjective 
homomorphism $\pi_1(O_n) \rightarrow
({\Bbb Z}/n{\Bbb Z}) \ast ({\Bbb Z}/n{\Bbb Z})$.
But $({\Bbb Z}/n{\Bbb Z}) \ast ({\Bbb Z}/n{\Bbb Z})$
has a free non-abelian subgroup of finite index, provided
$n > 2$. Therefore, $\pi_1(O_n)$ is large, and so
the same is true for $\pi_1(M(L,n))$, and hence
$\pi_1(O(L,n))$. $\square$

The following result allows us to consider some 3-orbifolds
with just one distinguished component of their
singular locus.

\noindent {\bf Theorem 3.2.} {\sl Let $O$ be a compact orientable
3-orbifold, and let $K$ be a knot in $O$ disjoint from
the singular set of $O$. Suppose that there is a surjective
homomorphism $\phi$ from $\pi_1(O)$ onto
some finite group $H$, so that $\phi(\langle [K] \rangle)$
has index at least $3$ in $H$. Then, provided
$n$ is sufficiently large, $\pi_1(O(K,n))$ is large.}

Here, $[K]$ is some element of $\pi_1(O)$ representing $K$,
and $\langle [K] \rangle$ is the subgroup generated
by $[K]$.
This is only defined up to conjugacy in $\pi_1(O)$, but
the index of $\phi(\langle [K] \rangle)$ in $H$ is nevertheless well-defined.

\noindent {\sl Proof of Theorem 3.2.} The kernel of $\phi$ corresponds to
a finite covering $\tilde O \rightarrow O$. Let $L$ be
the inverse image of $K$ in $\tilde O$. We then
have an induced finite covering $\tilde O(L,n) \rightarrow
O(K,n)$. The number
of components of $L$ is equal to the index of 
$\phi(\langle [K] \rangle)$ in $H$, which we are assuming is at least 3.
Hence, by Theorem 3.1, $\pi_1(\tilde O(L,n))$ is large
for all sufficiently large $n$, and the same is
therefore true for $\pi_1(O(K,n))$. $\square$

Such a homomorphism $\phi$ as in Theorem 3.2 very often exists.
For example, we shall show in Proposition 3.4 that this is always
the case if $O$ is a finite volume hyperbolic 3-orbifold. A variant of
Theorem 3.2 is as follows.

\noindent {\bf Theorem 3.3.} {\sl Let $O$ be a compact orientable 3-orbifold,
and let $K$ be a knot in $O$ disjoint from its
singular locus. Let $m$ be a positive integer such that
$\pi_1(O(K,m))$ admits a surjective homomorphism $\phi$ onto
a finite group $H$, with the property that $\phi(\langle [K] \rangle)$ 
has index at least $3m$ in $H$. Then, for all sufficiently
large $n$, $\pi_1(O(K,mn))$ is large.}

\noindent {\sl Proof.} This is very similar to the
proof of Theorem 3.2. The kernel of $\phi$
corresponds to a finite regular covering map $\tilde O \rightarrow O(K,m)$.
The inverse image of $K$ is a link $L$ in
$\tilde O$, with at least 3 components, by our hypothesis
on the index of $\phi(\langle [K] \rangle)$. Let
$q$ be the order of any singularity along $L$,
which we set to $1$ if $L$ is disjoint from the
singular locus of $\tilde O$. Let $\tilde O_n$ be the
orbifold with the same underlying manifold as
$\tilde O$, with singular locus consisting of that
of $\tilde O$, but with a singularity of
order $qn$ along $L$ (rather than $q$). We then have
an induced covering map $\tilde O_n \rightarrow
O(K,mn)$. Since $\pi_1(\tilde O_n)$ is large
for all sufficiently large $n$, by Theorem 3.1,
the same is true of $\pi_1(O(K,mn))$ for all sufficiently
large $n$. $\square$

The next result allows us to find homomorphisms
$\phi$ as in Theorems 3.2 and 3.3 in most cases
of interest.

\noindent {\bf Proposition 3.4.} {\sl Let $G$
be a finitely generated, residually finite group
that is not virtually cyclic. Then, for all
$g \in G$ and any integer $N$, there is a
surjective homomorphism $\phi$ from $G$ onto
a finite group $H$ such that $[H:\phi(\langle g \rangle)]$
is at least $N$.}

We will need the following elementary lemma.

\noindent {\bf Lemma 3.5.} {\sl Let $G$ be a finitely generated,
residually finite group that is not cyclic.
Then, $G$ has a finite index characteristic
subgroup $K$ such that $G/K$ is not cyclic.}

\noindent {\sl Proof.} Suppose, on the contrary,
that, for every finite index characteristic subgroup $K$
of $G$, $G/K$ is cyclic. Then $G$ must be abelian. For otherwise,
there are elements $g_1$ and $g_2$ of $G$ such that
$[g_1, g_2] \not=e$. By the assumption that 
$G$ is residually finite, there is a finite
index normal subgroup of $K_1$ of $G$
not containing $[g_1, g_2]$. Let $K$ be
the intersection of the images of $K_1$ under
all automorphisms of $G$. Then, $K$ is
a finite index characteristic subgroup of
$G$ not containing $[g_1,g_2]$. But,
we are assuming that $G/K$ must be
cyclic, which implies that $g_1K$ and
$g_2K$ commute. Hence, $[g_1,g_2]K = K$,
and so $[g_1, g_2] \in K$, a contradiction.
Therefore, $G$ is a finitely generated
abelian group. We will suppose that it is
not cyclic and reach a contradiction.
If it finite, then $\{ e \}$ is a finite
index characteristic subgroup such that
$G / \{ e \}$ is not cyclic, which is a 
contradiction. If $G$ is infinite, then it
is either of the form ${\Bbb Z}^n$, for some
$n \geq 2$, or 
${\Bbb Z} \times {\Bbb Z}/m{\Bbb Z}
\times H$, for some integer $m \geq 2$ and 
some abelian group $H$. In
both cases, set $K$ to be the subgroup
generated by the $m^{\rm th}$ powers of $G$
(where $m=2$, say, in the former case)
to achieve a contradiction. $\square$

\noindent {\sl Proof of Proposition 3.4.} We prove this by induction on $N$.
It is trivially true for $N = 1$. Suppose therefore,
that $N$ is at least 2, and that the inductive
hypothesis holds true for $N - 1$. This
implies that there is a surjective homomorphism
$\phi$ from $G$ onto a finite group $H$
such that $[H:\phi(\langle g \rangle)]$
is at least $N - 1$. Then, $K_1$, the kernel of $\phi$,
is finitely generated,
residually finite, and not cyclic.
Therefore, by Lemma 3.5, there is a 
finite index characteristic subgroup $K$
of $K_1$ such that
$K_1/K$ is not cyclic.
Now, $K$ is a finite index normal subgroup of
$G$. Let $\psi \colon G \rightarrow G/K$ be
the quotient homomorphism. We claim that the index
of $\psi(\langle g \rangle)$ in $G/K$ is at least $N$.
This index is $[G: \langle g \rangle K]$, which
equals $[G:\langle g \rangle K_1] [\langle g \rangle K_1
: \langle g \rangle K]$. The first of these terms
is, by assumption, at least $N - 1$. It therefore
suffices to show that $[\langle g \rangle K_1
: \langle g \rangle K]$ is at least two.
If this is not the case, then $\langle g \rangle K_1 = 
\langle g \rangle K$. Taking intersections with $K_1$,
we then deduce that $K_1 = \langle g \rangle K \cap K_1$.
This implies that $K_1 / K$ is cyclic,
generated by $g^nK$, where $g^n$ is a generator
for $\langle g \rangle \cap K_1$.
This is a contradiction. $\square$

Note that when $G$ is the fundamental group of
a finite volume hyperbolic 3-orbifold, then it
satisfies the hypotheses of Proposition 3.4:
it is finitely generated, residually finite and
not virtually cyclic. Hence, Theorem 3.3 and
Proposition 3.4 have the following
corollary.

\noindent {\bf Theorem 3.6.} {\sl Let $O$ be a compact orientable
3-orbifold (with possibly empty singular locus), and let $K$ be a knot in $O$, disjoint from its
singular locus, such that $O - K$ has a finite volume
hyperbolic structure. Then, for infinitely many
values of $n$, $\pi_1(O(K,n))$ is large.}

\noindent {\sl Proof.} It is a well known consequence of the proof of 
Thurston's hyperbolic Dehn surgery theorem [15] that
for all sufficiently
large $m$, $O(K,m)$ is hyperbolic. Hence, by
Proposition 3.4, $\pi_1(O(K,m))$ admits a
surjective homomorphism $\phi$ onto a finite group $H$, 
such that $\phi(\langle [K] \rangle)$ has index at least $3m$ in
$H$. Now apply Theorem 3.3 to deduce that
$\pi_1(O(K,mn))$ is large for all sufficiently
large $n$. $\square$

From the proof of the theorem, we obtain information
about the set (${\cal L}$, say) of values of $n$ for which 
$\pi_1(O(K,n))$ is large. Specifically,
there is an integer $A$ and, for each integer
$m \geq A$, an integer $B(m)$, such that
${\cal L}$ contains
$$\{ mn : m \geq A, n \geq B(m) \}.$$

We now focus on a classical case: cyclic branched covers
over a knot in the 3-sphere.

\noindent {\bf Theorem 3.7.} {\sl Let $K$ be a
non-trivial knot in the 3-sphere, and let $m$ be any integer
more than two. Then, for all sufficiently large $n$, the $mn$-fold
cyclic cover of $S^3$ branched over $K$ has
large fundamental group.}

\noindent {\sl Proof.} For any positive integer $n$, let $S^3(K,n)$ denote the
orbifold with underlying manifold $S^3$ and with
a singularity of order $n$ along $K$. Then the 
$n$-fold cyclic cover of $S^3$ branched over $K$
is a finite-sheeted covering space of $S^3(K,n)$. Our aim is therefore
to show that $\pi_1(S^3(K,n))$ is large for suitable
values of $n$.

Suppose first that $K$ is a connected sum of
two non-trivial knots $K_1$ and $K_2$. Then $S^3(K,n)$
is an orbifold connected sum of $S^3(K_1,n)$
and $S^3(K_2,n)$. Now, $\pi_1(S^3(K_1,n))$
and $\pi_1(S^3(K_2,n))$ are quotients of
$\pi_1(S^3(K,n))$. Hence if one of these has large
fundamental group, then so does $\pi_1(S^3(K,n))$.
So, it suffices to consider the case where $K$
is prime. This implies that $S^3(K,n)$ is (orbifold)-irreducible.

If $K$ is a satellite knot, then there is
an essential torus in its complement. Since $K$
is prime, this remains (orbifold)-incompressible
in $S^3(K,n)$, provided $n > 1$. So, its
inverse image $T$ in $M$, the $n$-fold 
cyclic branched cover of $S^3$ over $K$,
is incompressible, provided $n > 1$. 
Now, it is a theorem of Kojima [8] and Luecke [12]
that if a compact orientable irreducible 3-manifold $M$
contains essential embedded tori $T$, then
either $\pi_1(M)$ is large or $M$ is
finitely covered by a torus bundle over the circle,
with $T$ lifting to fibres. We claim that the
latter possibility cannot arise.
This is because one component
of the complement of $T$ covers a non-trivial
knot exterior, and this component would be covered by
$T^2 \times I$. However, Theorem 10.5 of [6] implies that
the only orientable irreducible 3-manifolds that are finitely covered by $T^2 \times I$
are $T^2 \times I$ itself and the orientable twisted $I$-bundle over a
Klein bottle. Therefore, $\pi_1(S^3(K,n))$ is large
when $n > 1$.

If $K$ is a $(p,q)$-torus knot, then $S^3(K,n)$ is
an orbifold Seifert fibre space with base orbifold
that is topologically a sphere and has
three singularities, with orders $p$, $q$ and $n$.
When $n$ is more than $6$, $(1/p) + (1/q) + (1/n) < 1$,
and therefore this base orbifold is hyperbolic.
Its fundamental group is therefore large. But
the Seifert fibration induces a surjective
homomorphism from $\pi_1(S^3(K,n))$ onto the
fundamental group of this orbifold, and therefore
$\pi_1(S^3(K,n))$ is large, when $n>6$.

Thus, we may assume that $K$ is hyperbolic.
Now, when $K$ is not the figure-eight knot,
it is a consequence of the Orbifold Theorem 
([2], Corollary 1.26 of [4]) that $S^3(K,m)$ is hyperbolic, when $m \geq 3$.
So, by Proposition 3.4 and Theorem 3.3, $\pi_1(S^3(K,mn))$
is large, when $n$ is sufficiently large. When $K$ is the figure-eight knot,
$S^3(K,m)$ is hyperbolic whenever $m \geq 4$,
and so the theorem also holds in this case. However,
$S^3(K,3)$ is Euclidean. Its fundamental group
is therefore residually finite and not virtually cyclic, and therefore
Proposition 3.4 and Theorem 3.3 combine 
to prove the theorem here also. $\square$

It is natural to speculate whether Theorems 3.6 and 3.7 can be strengthened.
Is it the case that when $K$ is a non-trivial knot in the
3-sphere, $\pi_1(S^3(K,n))$ is large for all but finitely
many values of $n$? This remains an interesting unsolved problem.
It suffices to consider the case where $n$ is prime,
but this is the main situation where the arguments in this
paper do not apply.

\vskip 18pt
\centerline{\caps References}
\vskip 6pt

\item{1.} {\caps B. Baumslag, S. Pride,} 
{\sl Groups with two more generators than relators,}
J. London Math. Soc. (2) 17 (1978) 425--426.

\item{2.} {\caps M. Boileau, J. Porti},
{\sl Geometrization of 3-orbifolds of cyclic type,}
Astérisque 272, (2001).

\item{3.} {\caps J. Button},
{\sl Strong Tits alternatives for compact 3-manifolds with boundary,}
J. Pure Appl. Algebra 191 (2004) 89--98.

\item{4.} {\caps D. Cooper, C. Hodgson, S. Kerckhoff},
{\sl Three-dimensional orbifolds and cone-manifolds.} MSJ Memoirs, 5. 
Mathematical Society of Japan, Tokyo (2000).

\item{5.} {\caps D. Cooper, D. Long, A. Reid,} 
{\sl Essential closed surfaces in bounded $3$-manifolds,} 
J. Amer. Math. Soc. 10 (1997) 553--563.

\item{6.} {\caps J. Hempel}, {\sl 3-Manifolds}, Ann. of Math. Studies 86, 
Princeton University Press (1976)

\item{7.} {\caps J. Howie}, {\sl Free subgroups in groups of 
small deficiency.} J. Group Theory 1 (1998) 95--112.

\item{8.} {\caps S. Kojima}, {\sl Finite covers of 3-manifolds
containing essential surface of Euler characteristic $=0$,}
Proc. Amer. Math. Soc. 101 (1987) 743--747.

\item{9.} {\caps M. Lackenby,} {\sl A characterisation of
large finitely presented groups}, 
J. Algebra. 287 (2005) 458--473.

\item{10.} {\caps D. Long,} {\sl Immersions and embeddings of
totally geodesic surfaces}, Bull. London Math. Soc. 19
(1987) 481--484.

\item{11.} {\caps D. Long, G. Niblo,}
{\sl Subgroup separability and $3$-manifold groups.}
Math. Z. 207 (1991) 209--215.

\item{12.} {\caps J. Luecke}, {\sl Finite covers of 3-manifolds containing
essential tori,} Trans. Amer. Math. Soc. 310 (1988) 381--391.

\item{13.} {\caps J. Morgan}, {\sl On Thurston's uniformization theorem 
for three-dimensional manifolds}, in The Smith
conjecture (New York, 1979) 37--125, Pure Appl. Math., 112, Academic Press.

\item{14.} {\caps J. Ratcliffe}, {\sl Euler characteristics of
3-manifold groups and discrete subgroups of ${\rm SL}(2, {\Bbb C})$},
J. Pure Appl. Algebra 44 (1987) 303--314.

\item{15.} {\caps W. Thurston}, {\sl The geometry and topology
of 3-manifolds}, Lecture Notes, Princeton, 1978.

\vskip 12pt
\+ Mathematical Institute, University of Oxford, \cr
\+ 24-29 St Giles', Oxford OX1 3LB, United Kingdom. \cr

\end